\documentclass[12pt]{article}       
\usepackage{a4}
\usepackage{amsmath,amssymb,amsthm} 
\usepackage[all]{xy}
\theoremstyle{plain}               
\newtheorem{defn}{Definition}[section]

\newtheorem{proposition}[defn]{Proposition}
\newtheorem{corollary}[defn]{Corollary}
\newtheorem{thm}[defn]{Theorem}
\theoremstyle{definition}          

\newtheorem{example}[defn]{Example}

\begin{document}
\title{The planar Tree Lagrange Inversion Formula}
\author{Lothar Gerritzen}
\date{29.09.03}
\maketitle

\textbf{Introduction}

\medskip

  A planar tree power series over a field $K$ is a formal expression
$$\sum c_T \cdot T$$ where the sum is extended over all isomorphism classes of finite planar reduced rooted trees
$T$ and where the coefficients $c_T$ are in $K$. Mulitplications of these power series is induced by planar grafting of trees
and turns the K-vectorspace $K\{x\}_\infty$ of those power series into an algebra, see [G].\\
If $f \in K \{x\}_\infty$  there is a unique $g(x) \in K \{x\}_\infty$ of order $> 0$ such that
$$ g(x) = x \cdot f(g(x))$$ where $f(g(x))$ is obtained by substituting $g(x)$ for $x$ in $f(x).$ Formulas for the coefficients
of $g$ in terms of the coefficients of $f$ are obtained by the use of the planar tree Lukaciewicz language.
This result generalizes the classical Lagrange inversion formula, see [C],[R],[Sch].

\section
{Lukaciewicz  languages}

Let $Y$ be a $\mathbb{N}$-graded set. It is given as a disjoint union of subsets $Y_i$ for $i \in \mathbb{ N}$.
Assume that $Y_0$ is not empty.\\
Let $\mathbb{W}(Y)$ be the monoid of words over $Y$ and let
$$\delta \colon \mathbb{W}(Y) \rightarrow \mathbb{N}$$
be the monoid morphism from $\mathbb{W}(Y)$ into additive monoid $\mathbb{N}$ of natural numbers for which
$$\delta(y) = i -1$$
whenever $y \in Y_i.$

\medskip

\begin{defn}
\textrm{Luk} $(Y)$ is the set of all words $w \in \mathbb{W}(Y)$ such that $\delta (w) = -1$ and such that if $w'$
 is a proper left factor of
$w$, then $\delta (w') \ge 0$.\\
 We call  \textrm{Luk}\ $(Y)$ the Lukaciewicz language over the graded set $Y$.
 \end{defn}

\medskip

In case\  $\# Y_i = 1$ for all $i$, then  \textrm{Luk}$(Y)$ is the classical Lukaciewicz language, for instance used in [C], (11.3).
We will show that the results in [C], (11.3) remain correct in the more general situation.

\medskip

Let $M(Y)$ be the submonoid of $\mathbb{W} (Y)$ generated by  \textrm{Luk}$(Y)$.

\medskip

\begin{proposition}

\smallskip

\begin{itemize}
\item[(i)]
$M(Y)$ is a freely generated by  \textrm{Luk}\ $(Y)$
\item [(ii)]
Let $w \in \mathbb{W}(Y)$ with $\delta(w) = - r, r \in \mathbb{N},$ then $w \in M(Y)$ if and only if $r > 0$ and if
$w'$ is a proper left factor of $w$, then $\delta(w')>-r$.
\end{itemize}
\end{proposition}

\begin{proof}
By definition  \textrm{Luk}$(Y)$ has the prefix property and therefore \textrm{Luk}$(Y)$ is a code which means that $M(Y)$
is freely generated by  \textrm{Luk}$(Y)$.
\end{proof}

\textbf{Example}
$$\bigcup^\infty_{i=1} x i$$
then $$\# Y_{0} = 1$$
$$Y_0 \in Y_0$$
$$ \mathbb{W} (Y +)\subseteq \textrm{Luk}$$
$$w \rightarrow w{r+1}_0$$
If $w_1,...,w_m \in \textrm{Luk}(Y)$ and $y_m \in Y_m, m \ge 0,$ then $v:=y_m \cdot w_1 \cdot ...\cdot w_m \in \textrm{Luk}
(Y)$\\
Obviously $$\delta(v)= \delta(y_m)+ \sum^m_{i=1}\delta (w_i)= $$
$$(m-1) + m(-1) = -1.$$ Let $v'$ be a proper left factor of $v, say   v = v' \cdot v". $
Then $v'= y_m w_1 \cdot ... \cdot w_k \cdot w'_{k+1}$ with $0 \le k < m,\  w'_{k+1}$ is a proper factor of $w_{k+1}$.\\
Then
$$\delta (v') = (m-1) + \sum^{k}_{i=1}  \delta (w'_{k+1})= $$

$$ =(m-k-1)+ \delta(w'_{k+1}) \ge m-k-1 \ge 0$$.

\begin{proposition}
Let $w \in\   \textrm{Luk}\ (Y)$ and let $w(1)$ be the first letter of $w$.\\
If $w(1) \in Y_m, m \ge 1,$
then there are unique $w_1,...,w_m \in\  \textrm{Luk}\ (Y)$ such that
$$ w= w(1) \cdot w_1 \cdot w_2,..., w_m$$
\end{proposition}

\begin{proof}
If $m=0$, then $w = w(1)$, because otherwise $w(1)$ would be a proper left factor of $w$ contained in $\textrm{Luk}\ (Y)$.\\
Let now $m > 0.$ Then $\delta(w') = \delta(w) - \delta(w(1))=(-1)-(m-1)=-m$.\\
There is a left factor $w'_1$ of $w'$ of smallest length relative to $Y$ such that
 $$\delta(w'_1)=-1$$
 Then $w'_1 \in\  \textrm{Luk}\ (Y)$. In a similar way $w'_2, ..., w'_m \in\  \textrm{Luk}\ (Y)$ are determined.

\end{proof}

 \medskip

 Next we define the height $h t (v)$ of an element $v \in \textrm{Luk}\ (Y)$ by induction on the length of $v$.\\
 If $v$ is a letter in $Y_0$, then $ h t (v):=0.$\\
 Let now $v \in \textrm{Luk}\ (Y)$ and assume that the length $l_Y(v)= n >1$. Then
 $$v= v(1) \cdot v'_1,...,v'_m$$
 if$$v(1) \in Y_m,\  m>0.$$
 and $l_Y(v'_i) < n $ for all $i$. By induction hypothesis we may assume that $ht(v'_i)$ is already defined for all $i$.\\
 Let $\displaystyle ht(v):=1+max^m_{i=1} ht (v'_i)$

\begin{defn}
$ht(v)$ is called the height of $v \in \textrm{Luk}\ (Y).$
\end{defn}
A classical result about Lukaciewicz languages  is described now.\\
Let  \textrm{Luk} $\mathbb{(N)}$ be the Lukaciewicz language relative to the decomposition of $\mathbb{N}$ into subsets
$ \mathbb{N}_i = \{i\}$ for $i \in \mathbb{N}$.
For $n \in \mathbb{N}$ we denote by $\lambda(n)$ the letter in $\mathbb{W} \in \mathbb{(N]}$ represented by $n$
 in order that no confusion occurs
between $\lambda(n) \cdot \lambda(m)$ and $\lambda(n \cdot m)$.\\
Let \textbf{PT} be the set of isomorphism classes of finite planar rooted trees, see [G 1] and \textbf{PT}$^\prime= \textbf{PT} \cup \{1\}$
where $1$ is the empty tree. The map $$w: \textbf{PT}^\prime \rightarrow \textrm{Luk} \mathbb{(N)}$$ is defined inductively with
respect to the number
$\# T^0$ of vertices of $T \in \textbf{PT}^\prime$.\\
If $T=1$, let $w(T):= \lambda(0) \in \textrm{Luk} \mathbb{(N)}$.\\ If $T = x$ is the tree in \textbf{PT} with a single vertex, then
$$w(x) = \lambda(1) \lambda(0) \in \textrm{Luk} \mathbb{(N)}.$$
If $\# T^0 >1$ and $\rho_T$ is the root of $T$, then the subgraph $T-{\rho_T}$ obtained by deleting the vertex
$\rho_T$ and all the edges in $T$ incident with $\rho _T$ is a forest. Any component $V_i$ of $T-{\rho_T}$ is a rooted
tree whose root is the vertex in $V_i$ incident with $\rho_T$ in $T$. Also the components $_1,...,V_m,\  m \ge 1,$ are
ordered and $\# V_i^0< \# T,$.\\
We may assume that $w(V_1),..., w(V_m)$ are already defined. Let $$w(T):= \lambda(m) \cdot w(V_1) \cdot...\cdot w(V_m).$$
Clearly $w(T) \in \textrm{Luk} \mathbb{(N)}$

\begin{proposition}
The map $w: \textbf{PT}^\prime \rightarrow \textrm{Luk} \mathbb{(N)}$ is bijective.\\
Moreover $ht(w(T)) = max_{b \in L(T)} (dist(\rho_T,b))$
where $dist(\rho_T,b)$ is the distance of the root $\rho_T$ of $T$ to a leaf $b$ in $T$ which is the length of a
smallest path in $T$ connecting $\rho_T$ and $b$.
\end{proposition}

\begin{proof}:1)
Let $S,T \in \textbf{PT}$ and assume that $w(S)=x(T).$\\
Then $$ar(S)=a(T)=m$$
If $$m=0$$
then $S=T$ is the empty tree.\\ If $m>0$, then
$$T=T_1 \cdot ... \cdot T_m$$
$$S=S_1 \cdot ... \cdot S_m$$
and\\
$\# T_i^0 < \# T^0, \#S_i^0 < \#S^0$. If we proceed by induction on $\#T^0 + \# S^0$ we see that $ T_i=S_i$
for all $i$. Thus $S=T$ and the map $w$ is proved to be injective.\\
2) Let now $v \in\  \textrm{Luk} \textbf{(PT)}^\prime$.
If $v \in Y_0$, then $v= w(1), 1$ empty tree. If the length of $v$ is greater than $1$, \\ than
$$v=v(1) \cdot v'_1 \cdot...\cdot v'_n$$
according to Proposition 1.4. Then $v(1)=\lambda (S),S \in \textbf{PT}$ and
$$ w(T_i)= v'_i$$ for $$T_1, ... , T_m \in \textbf{PT}.$$ Then $deg\ (S) = m$ and $w(*_S(T_1, ..., T_m))=v$ where
$*_S(T_1,...,T_m)$ is the planar grafting of $T_1,..., T_m$ over $S$.
\end{proof}

\begin{example}
$v= \lambda(3) \lambda(0) \lambda(1) \lambda(2) \lambda(0) \lambda(1) \lambda(0) \in \textrm{Luk} \mathbb{(N)}$\\
 and $w(T) = v$ for
\end{example}
$ T=$
\quad \quad \quad 
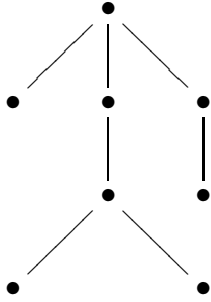
{\xymatrix{&\bullet \ar@{-}[dl] \ar@{-}[d]\ar@{-}[dr] &\\
\bullet & \bullet \ar@{-}[d] & \bullet\ar@{-}[d]\\
& \bullet \ar@{-}[dl] \ar@{-}[dr] &\bullet\\
\bullet && \bullet}

\section{The planar tree Lukaciewicz language}

Let \textbf{PRT} be the set of isomorphism classes of finite planar reduced rooted trees and \textbf{PRT}$^\prime= \textbf{PRT}
 \cup \{1\}$
where $1$ denotes the empty tree, see [G1]. Let\  $deg:  \textbf{PRT}^\prime \rightarrow \mathbb{N}$ be the degree map defined by
$deg(1)=0,\  deg(T) = \# L(T)$ where $L(T)$ is the set of leaves of $T$. We denote by Luk \textbf{(PRT}$^\prime$)
the Lukaciewicz
language relative to the grading by deg.\\
 In the following we will give an interpretation of the elements in Luk
 \textbf{(PRT}$^\prime)$ by right-sided decompositions of trees.
 For any tree $T \in \textbf{PRT} $ and any vertex $a$ of $T$, we denote by $T_a$ the closed subtree in $T$ generated by
 $a$.\\
 There is a unique edge $k$ in $T$ incident with $a$ such that $T_a$ is a connected component of $T-k$ and such that the
 root $\rho_T$ of $T$ is the root in $T_a.$
 This gives another characterization of $T_a.$\\
Let $T \in \textbf{PRT} $ and $x$ be the tree in \textbf{PRT} with the single vertex.

\begin{defn}
$T$ is called right-sided, if there is $T' \in \textbf{PRT}^\prime$ such that
 $$T=x \cdot T'$$
\end{defn}

 \bigskip
 The set of all right-sided trees in \textbf{PRT} is denoted by $x \cdot \textbf{PRT}!$\\
 Then
 $$ x=x \cdot 1,\  x^2=x \cdot x,\  x  \cdot x^2,\  x \cdot x^3,\  x \cdot(x^2 \cdot x),\  x \cdot(x \cdot x^2)\  \in\  x \cdot
 \textbf{PRT}!$$

\begin{defn}
$S$ is called completely right-sided in $T$, if\  $deg(S)>1$ and for any leaf $b \in L(S)$ the closed subtree $T_b$ in $T$
generated by $b$ is right-sided.\\

If $S$ is completely right-sided in $T$, then the forest
$$T-In(S)$$
obtained by recnoving all the inner vertices of $S$ from $T$ is a disjoint union of right-sided trees.

\medskip

Let $S,\  S'$ be completely right-sided open subtrees of $T$ and $S \subseteq S'$

\end{defn}

\begin{defn}
$S$ is strictly contained in $S'$, if relative to $T$ a leaf $b$ of $S$ which is also a leaf of $S'$ is already a leaf of $T$.
\end{defn}



\bigskip

Let $T \in x \cdot \textbf{PRT} $ and $r \ge 1.$\\
A right-sided open flag on $T$ of length $r$, is a sequence $$S=(S_1,\ S_2, ..., S_r)$$ with:
$$S_r=T$$
$S_i$ is a completely right-sided open subtree of $T$ for all $i<r\\ S_i$ is strictly contained in $S_{i+1}$ for all
$1 \ge i < r$.

\medskip

Denote by $\Lambda_r(T)$ the set of all proper right-sided open flags on $T$.\\
Let $\Lambda_r= \cup \Lambda_r(T)$ where the union is extended over all trees in $T$ for any $r \ge 1.$
Let $$\Lambda= \bigcup^\infty_{r=0} \Lambda_r$$ with
$$\Lambda_0= \{x\}$$
We are going to define a map $$w: \Lambda \longrightarrow \textrm{Luk}\textbf{(PRT}^\prime)$$
If $S \in \Omega_0,$ then $w(S):= \lambda(1):=$ letter in $\mathbb{W} \textbf{(PRT)}^\prime $ associated to the empty tree.\\
Let now
$$S \in \Lambda_r(T), r \ge 1, T \in x \cdot \textbf{PRT}$$
If $r=1,$ then $S=T$ and $w(T):= \lambda(T') \cdot \lambda(1)^m,$ if $T=x \cdot T',\  deg (T')= m$

\medskip

Let now $r \ge 2, S=(S_1,...,S_r).$ We proceed by induction on $deg(T).$\\
Consider the planar forest
$$T-In(S_1)=V_{0 \cup}V_1 \cup... \cup V_m$$
and denote the connected components of $T-In(S_1)$ by $V_0, V_1,...,V_m$ Then $V_0=x$ is the leaf of the left factor
of $T$. \\
Let
$$ S| V_i=(S_2 \cap V_i,...,S_r \cap V_i).$$
Then $$S_r \cap V_i= V_i$$ and $$S| V_i \in  \Lambda_{r-1}(V_i)$$ for all $i.$

Define\\
$$w(S):= \lambda(S'_1)\cdot w(S|V_1) w(S|V_2) \cdot...\cdot w(S|V_m)$$
where
$$S_1=x \cdot S'_1,\  S'_1 \in \textbf{PRT.}$$
As $w(S,V_i) \in $ \textrm{Luk} \textbf{PRT}$^\prime$\\
 and $deg(S_1')=m$ we get $w(S) \in\    \textrm{Luk} \textbf{(PRT}^\prime)$

\begin{thm}
The map $$ w:\Lambda \rightarrow \textrm{Luk} \textbf{(PRT}^\prime)$$
is bijective.
\end{thm}

\section{Right-sided decompositions}
Let $T \in \textbf{PRT}$ and $S$ be a non-empty subtree of $T$. The vertex in $S$ closest to the root of $T$ is defined to be the
root $\rho_S$ of $S$ which turns $S$ into a rooted tree.
\begin{defn}
$S$ is relatively open in $T$, if $S$ is an open subtree in the closure $T_{(\rho_s)}$ of $\rho s$ in $T$.
\end{defn}

It is easy to show that $S$ is relatively open in $T$ of and only if $ar_S(a)= ar_T(a)$ for all vertices $a$ of $S$ with
$ar_S(a)>0.$ There $ar_T(a)$ denotes the number of outgoing edges of the vertex $a$ in $T.$\\

Let $T$ be a right-sided tree in \textbf{PRT} of degree $ > 1, T=x T',$ and let $Q$ be a system of relatively open subtrees
$U$ of $T, T=x \cdot T'$

\begin{defn}
$Q$ is called right-sided decomposition of $T$, if
\begin{itemize}
\item[(i)]
Each $U \in Q$ is right-sided
\item [(ii)]
$\displaystyle \bigcup_{U \in Q} U=T$
\item [(iii)]
If $b$ is a leaf of $T$, then $\{b\} \in Q$ if and only if $b$ is different from the first leaf of $T.$
\item[(iv)]
If  $U \cap U' \not= \phi , U \not= U',$ then $U \cap U'=\{a\},$ a vertex of  $T$, and if a is not the root of  $U$, then
it is a leaf in $T$ and the root of  $U'.$
\end{itemize}
\end{defn}
Denote by $\Gamma(T)$ the set of all right-sided decompositions of $T.$ Recall that $\Lambda(T)$ denotes the set of all
proper open flags of $T.$ For $S=(S_1,...,S_r) \in \Lambda(T)$ let $Q(S)$ denote the systems of subtrees of $T$  which
contains $S_1,$ all the connected components of $S_i-In(S_{i-1})$ for $1<i \le r$ and $\{b \in L(T):b \not=$  first
leaf of $T \}:=L(T') \}$ , if   $T=x \cdot T'.$\\
Then
$$Q(S) \in \Gamma(T)$$

\begin{proposition}
The map $$S \longmapsto Q(S)$$
is a bijection
$$\Lambda(T) \rightarrow \Gamma(T)$$
for all $$T \in x \cdot \textbf{PRT.}$$
\end{proposition}

\textbf{Example}\quad \quad \quad \quad\\

$T = x \cdot (x \cdot x^2 \cdot x \cdot x^3)$

$
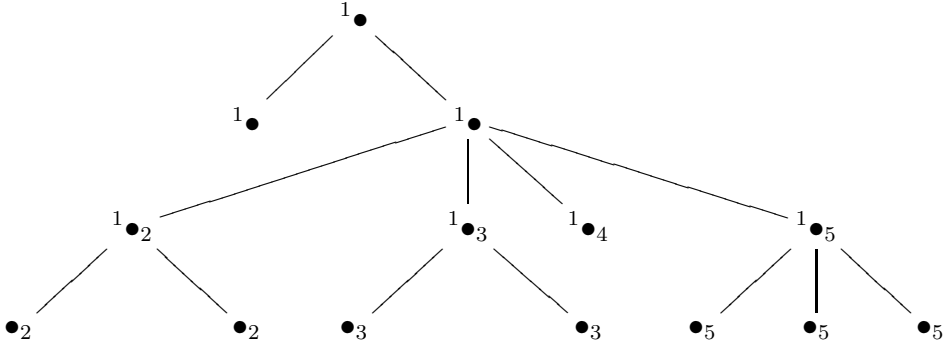
\xymatrix{&&&^1\bullet\ar@{-}[dr] \ar@{-}[dl] &&&&&\\ &&^1\bullet
&&^1\bullet \ar@{-}[drrr] \ar@{-}[d]
\ar@{-}[dlll]\ar@{-}[dr]&&&&\\ &^1\bullet_2\ar@{-}[dl]\ar@{-}[dr]
&&&^1 \bullet_3 \ar@{-}[dl] \ar@{-}[dr] & ^1\bullet_4&&
^1\bullet_5 \ar@{-}[dl] \ar@{-}[d] \ar@{-}[dr]&\\ \bullet_2
&&\bullet_2& \bullet_3 &&\bullet_3  &\bullet_5 &\bullet_5
&\bullet_5}$

\bigskip

$S=$ subgraph of $T'$ vertices indexed by $1$.
$S$ is open completly right-sided subtree of $T$ and $T-In(S)= V_2 \cup...\cup V_5$ where $V_i=$ subgraph of $T$ of vertices
indexed by $i$.

\section{Lagrange inversion formula for planar tree power series}
 Let $K$ be a field and $\textbf{A}=K \{\{x\}\}_\infty$ the K-algebra of planar tree power
series in $x$ over $K.$
\begin{proposition}
Let $$f \in \textbf{A},\   ord(f)=0$$
Then there is a unique power series
$$1/f \in \textbf{A}$$
such that $$f \cdot(1/f)=1.$$
Then $ord(x \cdot(1/(f))=1$ and there is a unique $g \in \textbf{A}$ such that
$$ord(g)=1$$
$$ \varphi_g (S x 1/f)=x$$
\end{proposition}
\begin{proof}1)
The coefficients $\gamma(T)$ of $1/f$ satisfy the following system of equations:
$$c_1(f) \cdot \delta(1) = 1$$
$$c_x(f) \cdot \gamma(1) + c_1(f) \cdot \gamma(x)=0$$
If $$ar(T)=m \ge 2,  \cdot T=T_1 \cdot T_2 \cdot ...\cdot T_m$$
Then
$$c_T(F) \gamma(1) + c_1(f) \gamma(T)=0$$
if
$$ m>2.$$
If $$m=2$$
Then
$$c_T(f) \gamma(1) +c_1(f) \cdot \gamma(T) + c_{T_1}(f) \gamma(T_2)=0$$
There is a unique solution of this system of equation and $\sum \gamma(T) \cdot T$ is equal to $1/f.$\\
2) $\varphi_{(x \cdot 1/f)}$ is an automorphism of $\textbf{A}$ as linear term of $x \cdot 1/f \not= 0.$
Thus the inverse of $\varphi(x \cdot 1/(f)$ exists and is equal to $\varphi_g$ for some $g \in \textbf{A}$
Then
$\varphi_g$ for some $g \in \textbf{A}$
Then
$$\varphi_g\circ \varphi_{(x\cdot 1/f)}=\varphi_x$$
and
$$ \varphi_{g(x \cdot 1/f)} = x.$$
\end{proof}

\begin{proposition}
Let $f \in \textbf{A}.$ Then there is a unique $g \in \textbf{A}$ such that $g(x)=x \cdot f(g(x)).$ It follows
that $ord(g) \ge 1.$
\end{proposition}

\begin{proof}1)
Let $a(t)=c_T(f)$
Now there is a unique system\  $b(T)$\  of coefficients such that
$$b(T)=0\ \textrm{if}\  T \not= x \cdot T'$$
$$b(x \cdot T')=\sum_{S' \in \Omega(T')}a(S') b(T' - In(S'))$$
Then
$$ g= \sum b(T) \cdot T$$
satisfies the equation of the Proposition.
2)Let
$b(T)= c_T(g) \cdot$
Then
$$ c_T(x \cdot f(g(x))= 0$$ if $$T \not= x \cdot T'$$
and
 $$ c_{x \cdot T'} (x \cdot f(g(x))= c_{T'} (f(g(x))$$
 $$=\sum_{S' \in \Omega(T')} b(T'-In(S'))$$
 $$=\sum_{S \in \Omega (T)\atop S \not= \delta_T} a(S') b(T'-In(S))$$
 Let
  $$a(T) = c_T(f).$$

\end{proof}

\begin{corollary}:
$$g(x)= x \cdot f(g(x))$$
\end{corollary}

\begin{proof}
$$\varphi_g(x \cdot 1/f)= g(x) \cdot (1/f) (g(x))$$
but
$$(1/f(g(x)= 1/f(g(x)$$
because
$$ord(f(g(x)=0.$$
\end{proof}

\begin{thm}
Let $f, g$ as in the Proposition above and
$$a(T)=c_T(f),\  b(T)= c_T(g).$$
Then
$$ b(x)= a(1)$$
and if $$T \in x \cdot \textbf{PRT},\  deg(T) > 1,$$
then
$$ b(T)= \sum_{Q \in \Gamma(T)}a(Q)$$
where
$$ a(Q):= \prod_{U \in Q \atop U=x \cdot U'}a(U').$$
\end{thm}

\begin{proof}
We proceed by induction on $n=\  deg(T).$ If $n=1,$ it is clear.\\
 If $n>1,$ then
$$b(T)= \sum a(S) \cdot b(T-In(S))$$
where the summation is extended over all proper completely right-sided open subtrees $S$ in $T$.
By induction hypothesis
$$b(T-In(S))=\prod_{i=1}^m b(V_i)$$
if $$T-In(S)= V_1 \cup ... \cup V_m$$
 and
$$b(V_i)= \sum_{Q \in \Gamma(V_i)}a(Q)$$
The map $Q(T) \rightarrow Q(V_1) \times...\times Q(V_m)$ which maps $S$ onto $(S(V_1,...,S L (V_m))$ is a bijection.
This proves the formula.
\end{proof}


\begin{thebibliography}{ccccc}
\bibitem[C]{C} Cori, R,:\  Words and Trees, in M. Lothaire, Combinatorics on Words, Cambridge 1981, Chap 11, 213 - 227
Manuskript 2003
\bibitem[R]{R} Raney, G. N. 1960:\ Functional composition patterns and power series reversion, Trans. Am. Math. Soc., 94, 441-451
\bibitem[Sch] {Sch} Sch\"utzenberger, M. P. 1971:\ Le th\`eor\`em de Lagrange selon Raney, in Logiques et Automates, Publications
Inst. Rech. Informatique et Automatique, Rocquencourt, France.
\bibitem[G 1] {G 1} Gerritzen, L.:\ Planar rooted trees and non-associative exponential series, Advances in Applied Mathematics 33
(2004) 342-365
\bibitem[G 2] {G 2} Gerritzen, L.:\ Automorphisms of the planar tree power series algebra and the non-associative logarithm,
Serdica Math. J. 30 (2004), 135-158
\bibitem [G 3] {G 3} Gerritzen, L.:\ Planar Shuffle Product, co-Addition and the non-associative Exponential,
ArXiv:math.RA/0502378
\end{thebibliography}
\end{document}